\documentclass[12pt,a4paper]{article}
\usepackage{latexsym, amssymb, amsthm}
\begin{document}
\theoremstyle{plain}
\newtheorem{theorem}{Theorem}[section]

\newtheorem{lemma}[theorem]{Lemma}
\newtheorem*{Cor}{Corollary}
\newtheorem{corollary}[theorem]{Corollary}
\mathsurround 2pt

\title{Products of locally cyclic groups}
\author{Bernhard Amberg and Yaroslav Sysak\thanks
{The second author likes to thank the Institute of Mathematics of the University of Mainz for its excellent hospitality during the preparation of this paper.}}
\renewcommand{\thefootnote}{\fnsymbol{footnote}}
\setlength{\footnotesep}{.1in} 
\date{}

\bibliographystyle{plain}

\maketitle
\begin{abstract}
We consider groups of the form G = AB with two locally cyclic subgroups A and B. The structure of these groups is determined in the cases when A and B are both periodic or when one of them is periodic and the other is not. Together with a previous study of the case where A and B are torsion-free, this gives a complete classification of all groups that are the product of two locally cyclic subgroups. As an application, it is shown that the Pr\"ufer rank of a periodic group product of two locally cyclic subgroups does not exceed 3. It is also proved that the product of a finite number of pairwise permutable periodic locally cyclic subgroups is a locally supersoluble group. This generalizes a well-known theorem of B. Huppert for finite groups.\end{abstract}

\section{Introduction}

Let the group $G=AB$ be the product of two subgroups $A$ and $B$, i.e. $G=\{ab\mid a\in A, b\in B\}$. It was proved by N. It\^ {o} that the group $G$ is metabelian if the subgroups $A$ and $B$ are abelian (see \cite[Theorem 2.1.1]{AFG}). This result laid the foundation for a systematic study of groups of the form $G = AB$ with various conditions on the subgroups $A$ and $B$. In particular, it follows directly from It\^o's result that every periodic group $G=AB$ with abelian subgroups $A$ and $B$ is locally finite. It is also well known that the group $G=AB$ with cyclic subgroups $A$ and $B$ is supersoluble and abelian-by-finite (\cite[Lemma 7.4.6]{AFG}). Furthermore, a detailed description of the structure of a group $G=AB$ with torsion-free locally cyclic subgroups $A$ and $B$ was obtained by the second author in \cite{S_86}. 

The aim of this paper is to describe the structure of groups that are products of two locally cyclic subgroups in the periodic and in the mixed case. Altogether this gives a complete answer to Question 15 in the book \cite {AFG}. 

\begin{theorem}\label{0} Let the periodic group $G=AB$ be the product of two locally cyclic subgroups $A$ and $B$. Then $G$ contains uniquely determined locally cyclic normal subgroups $S$ and $T$ and a locally nilpotent subgroup $H=A^*B^*$ with $A^*\le A$ and $B^*\le B$ such that 
$$G=(S\times T)\rtimes H=(S\times A^*)(T\times B^*)$$ where $\pi(S)\cap \pi(A^*)=\pi(T)\cap \pi(B^*)= \emptyset$, $S=[S,B^*]$ and $T=[T,A^*]$.\end{theorem}

We recall that a group $G$ has finite Pr\"ufer rank $r=r(G)$ (special rank in the sense of Mal'cev in Russian terminology) if each finitely generated subgroup of $G$ can be generated by $r$ elements and $r$ is the least positive integer with this property. Clearly a group is of rank $1$ if and only if it is locally cyclic. It is also known that every finite $p$-group of the form $G=AB$ with cyclic subgroups $A$ and $B$ has rank at most $2$ for $p$ odd (\cite{Hup}, Satz 8) and at most $3$ for $p=2$ (\cite{Ja}, Theorem 5.1). The following consequence of Theorem \ref{0} gives an exact upper bound for the Pr\"ufer rank of a product $G=AB$ of two periodic locally cyclic subgroups $A$ and $B$. 

\begin{corollary}\label{01} If $G=AB$ is a periodic group with locally cyclic subgroups $A$ and $B$, then the Pr\"ufer rank of $G$ does not exceed $3$.\end{corollary}

It should be noted that this result is also new for arbitrary finite groups of the form  $G=AB$ with cyclic subgroups $A$ and $B$. 

Finally, the following theorem extends a well-known result of B. Huppert on the supersolubility of finite groups which are products of pairwise permutable cyclic subgroups (see \cite{Hup}, Satz 34 or \cite{Hu}, Satz VI.10.3). This gives, in particular, an affirmative answer to Question 1 of the article \cite{AS_16}.

\begin{theorem}\label{12} Let the group $G=A_1A_2...A_n$ be the product of finitely many pairwise permutable periodic locally cyclic subgroups $A_1,...,A_n$. Then $G$ is a periodic locally supersoluble group.\end{theorem}

\section{Preliminaries} 

In what follows $G=AB$ is a group with locally cyclic subgroups $A$ and $B$. 

\begin{lemma}\label{1} If $G=AB$ is an infinite $p$-group, then up to a  permutation of the factors $A$ and $B$ the subgroup $A$ is quasicyclic and one of the following statements holds: \begin{itemize}
\item[1)] $G=A\times B$ with $B$ cyclic or quasicyclic;
\item[2)] $p=2$ and $G=A\rtimes \langle b\rangle$ with $a^b=a^{-1}$ for every $a\in A$;
\item[3)] $p=2$ and $G=A\langle b\rangle$ with $b^{2^n}=1$ for some $n>1$, $b^{2^{n-1}}\in A$ and $a^b=a^{-1}$ for every $a\in A$.\end{itemize}\end{lemma}

\proof Clearly without loss of generality we may assume that the subgroup $A$ is infinite and so quasicyclic. Then the subgroup $B$ is either cyclic or quasicyclic. Since in the latter case the group $G$ is abelian by \cite{AFG},  Lemma 7.4.4, the subgroup $A$ is complemented in $G$ and hence statement 1) holds. 

Let the group $G$ be non-abelian. Then the subgroup $B$ is cyclic and so $A$ as a quasicyclic $p$-subgroup of finite index in $G$ must be normal and non-central in $G$. In particular, $B$ induces on $A$ a non-trivial cyclic $p$-group of automorphisms. On the other hand, since quasicyclic $p$-groups have no automorphisms of order $p>2$, it follows that $p=2$. But then $B=\langle b\rangle$ with $b^{2^n}=1$ for some $n\ge 1$ and $b$ induces on $A$ an automorphism of order $2$ that inverts the elements of $A$. In particular, if $A\cap B=1$, we obtain statement 2). In the second case, $A\cap B=\langle b^{2^{n-1}}\rangle$ and hence statement 3) holds, as claimed.\qed

\begin{corollary}\label{2} If $G=AB$ is a $p$-group and $C, D$ are subgroups of $A$ and $B$, respectively, then $CD=DC$ and so $CD$ is a subgroup of $G$.\end{corollary}

\proof This is known if $G$ is finite (see \cite{Hup}, Satz 3), and follows from Lemma \ref{1} in the general case.\qed

\begin{lemma}\label{3} If the group $G=AB$ is periodic and $H$ is a finite subgroup of $G$, then $H$ is contained in a finite subgroup $E$ of $G$ such that $E=(A\cap E)(B\cap E)$. In particular, the group $G$ is locally supersoluble.\end{lemma} 

\proof Since the group $H$ is finite, there exist finite subsets $C$ of $A$ and $D$ of $B$ such that $H$ is contained in the set $CD$. Then the subgroups $A_0=\langle C\rangle$,  $B_0=\langle D\rangle$ and $\langle C,D\rangle$ are finite, because the group $G$ is locally finite. 
Furthermore, it follows from \cite{AFG}, Lemma 1.2.3, that the normalizer $N_G(\langle A_0, B_0\rangle)$ contains a finite subgroup $E$ such that $\langle A_0, B_0\rangle\le E=(A\cap E)(B\cap E)$. Since the subgroup $E$ is supersoluble by \cite[Lemma 7.4.6]{AFG} and $H\subseteq CD\subseteq\langle C,D\rangle = \langle A_0, B_0\rangle$, the lemma is proved.\qed\medskip

As a direct consequence of this lemma, we have 

\begin{corollary}\label{4} If the group $G=AB$ is periodic, then there exists an ascending series of  finite subgroups $1=G_0<G_1<\dots<G_n<\dots G$ such that $G_n=(A\cap G_n)(B\cap G_n)$ for each $n>0$ and $G=\bigcup_{n=1}^\infty G_n$. \end{corollary}

If $G$ is a periodic group and $\pi$ is a set of primes, then a subgroup $H$ of $G$ is called a $\pi$-subgroup, provided that all prime divisors of the order of any element of $H$ are contained in $\pi$.  
By a Sylow $\pi$-subgroup of $G$ we mean simply a maximal $\pi$-subgroup $G_{\pi}$ of $G$ which will be denoted by $G_p$ if $\pi=\{p\}$.

\begin{lemma}\label{5} Let $G=AB$ be a periodic group and $\pi$ a set of primes. Then the following statements hold.\begin{itemize}
\item[1)] If $A_{\pi}$ and $B_{\pi}$ are Sylow $\pi$-subgroups of $A$ and $B$, respectively, then $G_{\pi}=A_{\pi}B_{\pi}$ is a Sylow $\pi$-subgroup of $G$ and $$N_G(G_{\pi})=N_A(G_{\pi})N_B(G_{\pi}).$$
\item[2)] If $p, q$ are primes with $p>q$, then a Sylow $p$-subgroup $G_p$ is normalized by a Sylow $q$-subgroup $G_q$. In particular, $G_{\{p, q\}}=G_pG_q$ for any primes $p$ and $q$.\end{itemize}\end{lemma}

\proof 1) It follows from \cite{AFG}, Lemma 1.3.2, that in the notation of Corollary \ref{4} for each $n\ge 1$ the set $(A_{\pi}\cap G_n)(B_{\pi}\cap G_n)$ is a Hall ${\pi}$-subgroup of $G_n$. Therefore $G_{\pi}=\bigcup_{n=1}^\infty (A_{\pi}\cap G_n)(B_{\pi}\cap G_n)$ is a Sylow $\pi$-subgroup of $G$. Since  $A_{\pi}=\bigcup_{n=1}^\infty (A_{\pi}\cap G_n)$ and $B_{\pi}=\bigcup_{n=1}^\infty (B_{\pi}\cap G_n)$, this implies $G_{\pi}=A_{\pi}B_{\pi}$. In addition, applying \cite{AFG}, Lemma 1.2.2, we have  $N_G(G_{\pi})= N_A(G_{\pi})N_B(G_{\pi})$.

2) If ${\pi}=\{p, q\}$ for some primes $p>q$, then $A_{\{p,q\}}=A_p\times A_q$, $B_{\{p,q\}}=B_p\times B_q$ and $G_{\{p,q\}}=(A_p\times A_q)(B_p\times B_q)$. As $G$ and so its subgroup $G_{\{p,q\}}$ are locally supersoluble by Lemma \ref{3}, the Sylow $p$-subgroup $G_p=A_pB_p$ is normal in $G_{\{p,q\}}$. Therefore, $G_{\{p,q\}}=(A_p\times A_q)(B_p\times B_q)=G_pA_qB_q=G_pG_q$, as claimed.\qed\medskip

A Sylow basis of a periodic group $G$ is defined to be a complete set ${\bf S} =\{G_p\}$ of Sylow $p$-subgroups of $G$, one for each prime $p$, such that $G_pG_q=G_qG_p$ for all pairs $p, q$ of primes, and $G_{\pi}=\langle G_p \mid p\in{\pi}\rangle$ is a Sylow $\pi$-subgroup of $G$ for each set $\pi$ of primes. 
As is well known (see \cite{Har},Lemma 2.1), every countable periodic locally soluble group possesses Sylow bases. The basis normalizer $N_G({\bf S})$ of a Sylow basis $\bf S$ of $G$ is by definition the intersection $N_G({\bf S})=\bigcap_p N_G(G_p)$ of the normalizers $N_G(G_p)$ of the Sylow $p$-subgroups $G_p$ of $\bf S$ for all $p$.

\begin{lemma}\label{6} Let the group $G=AB$ be periodic and $G_p=A_pB_p$ for each prime $p$. Then ${\bf S}=\{G_p\}$ is a Sylow basis of $G$. Moreover, if $A^*=\bigcap_p N_A(G_p)$ and $B^*=\bigcap_p N_B(G_p)$, then $N_G({\bf S})=A^*B^*$.\end{lemma} 

\proof Indeed, by Lemma \ref{5} the set ${\bf S}=\{G_p\}$ forms a Sylow basis of $G$ and $N_G(G_p)= N_A(G_p)N_B(G_p)$ for every $p$ by \cite{AFG}, Lemma 1.2.2. Therefore, $N_G({\bf S})=\bigcap_p N_G(G_p)=\bigcap_p N_A(G_p)N_B(G_p)$ and  it is easy to check that $\bigcap_p N_A(G_p)N_B(G_p)=(\bigcap_p N_A(G_p))(\bigcap_p N_B(G_p))=A^*B^*$ (see \cite{AFG}, Lemma 1.1.2). Therefore $N_G({\bf S})=A^*B^*$.\qed\medskip

The following lemma is a direct consequence of a well-known result of L. Kovacs (see \cite{Ko}, Theorem 2).

\begin{lemma}\label{7} Let $G$ be a finite soluble group, $\pi$ a set of primes and $H$ a Hall $\pi$-subgroup of $G$. If for each $p\in\pi$ the Pr\"ufer rank of a Sylow $p$-subgroup of $G$ does not exceed $r$, then $H$ is a subgroup of rank at most $r+1$.\end{lemma} 

\proof Indeed, it is obvious that if $K$ is a subgroup of $H$, then every Sylow subgroup of $K$ is generated by $r$ elements. Therefore, $K$ can be generated by $r+1$ elements by the result of Kovacs cited above. Thus every subgroup of $H$ is generated by $r+1$ elements and so $H$ has rank at most $r+1$.\qed

\section{Proof of Theorem 1.1} 

First of all, it follows from Lemma \ref{1} that for each prime $p$ every Sylow $p$-subgroup of $G=AB$  satisfies the minimal condition for subgroups. Therefore, $G$ satisfies the minimal condition for $p$-subgroups for all primes $p$. Since the group $G$ is metabelian by Ito's theorem, the locally nilpotent residual $R$ of $G$ is contained in its derived subgroup $G'$ and so is abelian. It was proved by B. Hartley in \cite{Har}, Theorem 1, that in this case $G=R\rtimes H$, where $H$ is any basis normalizer of $G$. In particular, by Lemma \ref{6} we can take $H=A^*B^*$. 

It is easy to see that the subgroup $H$ is locally nilpotent and contains the center $Z(G)$ of $G$. Furthermore, $G'=R\times H'$ and so $H'$ is a normal subgroup of $G$. Since $R$ is abelian and $N_G(H)= N_R(H)\times H$, it follows that $N_R(H)\le Z(G)\le H$. Therefore, $N_R(H)=1$ and hence $H=N_G(H)$.  We show now that the subgroup $H=A^*B^*$ commutes with both subgroups $A$ and $B$.

Indeed, put $S=R\cap\langle A,B^*\rangle$ and $T=R\cap\langle A^*,B\rangle$. It is clear that $S$ and $T$ are normal subgroups of $G$, $\langle A,B^*\rangle=S\rtimes H$ and $\langle A^*,B\rangle=T\rtimes H$. On the other hand, as $G=AB$, we have also  $\langle A,B^*\rangle=AB_1$ and $\langle A^*,B\rangle=A_1B$ for some subgroups $A_1$ and $B_1$ such that $A^*\le A_1\le A$ and $B^*\le B_1\le B$. From here we deduce $AB_1\cap A_1B=A_1B_1=(S\cap T)\rtimes H$. 
Moreover, passing to the factor group $G/H'$, we may restrict ourselves to the case when the subgroup $H=A^*B^*$ is abelian. Then the subgroups $A^*$ and $B^*$ centralize $S$ and $T$, respectively, and so the subgroup $H$ centralizes the intersection $S\cap T$. Since $H=N_G(H)$, this implies $S\cap T=1$. Thus  $A_1B_1=H=A^*B^*$ and hence $\langle A,B^*\rangle=AH=AB^*$ and $\langle A^*,B\rangle=BH=A^*B$, as asserted. 

Further, taking into account the equalities $AB^*=S\rtimes H$ and $H=A^*B^*$, we conclude that the subgroup $A^*$ centralizes $S$, because $[A^*,S]\le H'\cap S=1$. Since in this case the normalizer $N_S(B^*)$ is contained in $N_G(H)=H$, we have $N_S(B^*)=1$. Therefore, every element $b\in B^*$ induces on $S$ an automorphism leaving only the identity element fixed. But then every element of $S$ can be written in the form $b^{-1}s^{-1}bs$ with $s\in S$ and hence $S=[B^*,S]$. Similarly, using the equality $A^*B=T\rtimes H$, we derive $T=[A^*,T]$.

Finally, we put $A_0=A\cap BS$ and $B_0=AT\cap B$. Clearly from the equalities $G=AB$, $AB^*=S\rtimes H$ and $A^*B=T\rtimes H$ it follows that $G=S\times A^*B=T\times AB^*$, $A=A^*\times A_0$ and $B=B^*\times B_0$.   Therefore, $S\rtimes B=A_0B$ and $T\rtimes A=AB_0$. Furthermore, if $S_p$ is a Sylow $p$-subgroup of $S$, then $S_p\rtimes B=(A_0\cap S_pB)B$ and in particular $S_p\ne 1$ if and only if $A_0\cap S_pB\ne 1$. Since the subgroup $A_0$ is locally cyclic, this implies that also $S$ is locally cyclic and $\pi(S)=\pi(A_0)$. Moreover, as $A^*$ and $A_0$ are subgroups of the locally cyclic subgroup $A$, it also follows that $\pi(A^*)\cap\pi(A_0)= \emptyset$. Similarly, using the equality $T\rtimes A=AB_0$, we obtain $\pi(T)=\pi(B_0)$ and $\pi(B^*)\cap\pi(B_0)=\emptyset$.\qed\medskip

{\bf Proof of Corollary 1.2}. By Corollary \ref{4}, we can restrict ourselves to the case in which the group $G=AB$ is finite. By Theorem \ref{0} $G$ contains cyclic normal subgroups $S$ and $T$ and a nilpotent subgroup $H=A^*B^*$ with $A^*\le A$ and $B^*\le B$ such that 
$$G=(S\times T)\rtimes H=(S\times A^*)(T\times B^*),$$ where $\pi(S)\cap \pi(A^*)=\pi(T)\cap \pi(B^*)= \emptyset$, $S=[S,B^*]$ and $T=[T,A^*]$. In particular, if for some prime $p$ the subgroup $H$ contains a non-cyclic Sylow $p$-subgroup $P$, then $S$ and $T$ are $p'$-subgroups of $G$. 

Since $P=A_pB_p$ with $A_p=A\cap P$ and $B_p=B\cap P$, both subgroups $A_p$ and $B_p$ are non-trivial and so $p\notin\pi(S)\cup\pi(T)$. Therefore, if $G_p$ is a non-cyclic Sylow $p$-subgroup of $G$, $S_p=G_p\cap S$ and $T_p=G_p\cap T$, then up to conjugation $G_p$ coincides with one of the following subgroups of $G$: $P=A_pB_p$,  $A_p\rtimes T_p$, $B_p\rtimes S_p$ and $S_p\times T_p$. In particular, the Sylow $p$-subgroups of $G$ have rank at most $2$ for $p>2$ (see \cite{Hu}, Satz III.11.5) and at most $3$ for $p=2$ (see \cite{Ja}, Theorem 5.1). We show now that $G$ is in fact a group of rank at most $3$. 

Indeed, suppose the contrary and let the group $G$ contain a subgroup $K$ whose minimal number of generators $d(K)$ is at least $4$. Since the Sylow subgroups of odd orders in $K$ have rank at most $2$ by what was noted above, each Sylow $2$-subgroup $Q$ of $K$ must have rank $3$ by Lemma \ref{7}. Furthermore, there exists a Sylow $2$-subgroup $P=A_2B_2$ of $G$ such that $Q=K\cap P$. As the group $G$ is metabelian, the derived subgroup $P'$ is normal in $G$ and so the subgroup $N=Q\cap P'$ is normal in $K$. It is easy to see that $N\ne 1$, because otherwise the subgroup $Q$ is embedded in the factor group $P/P'$ whose rank is equal $2$. Furthermore, the factor group $Q/N=Q/Q\cap P'$ is isomorphic to the factor group $QP'/P'\le P/P'$ and so the rank of $Q/N$ does not exceed $2$. Thus the factor group $K/N$ has rank at most $3$ by Lemma \ref{7}. In particular, $d(K/N)<d(K)=4$ and hence $N$ is not contained in the Frattini subgroup $\Phi(K)$ of $K$, because otherwise $d(K/N)=d(K)$. Therefore, passing to the factor group $K/\Phi(K)$, we may assume that $\Phi(K)=1$. Then the normal subgroup $N$ is complemented in $K$ and so in $Q$ by \cite{Hu}, Hilfsatz 3.2.b). 

Let $L$ be a complement to $N$ in $K$ and $M=Q\cap L$. Then $M$ is a Sylow $2$-subgroup of $L$ and $Q=M\times N$, so that $d(M)\le 2$ and the subgroup $Q$ is abelian. Therefore, the subgroup $L$ is three-generated by Lemma \ref{7} and $N$ is a central subgroup of $K$. In particular, if $d(M)=1$ and $N=\langle u, v\rangle$,  then generators $a, b, c$ of $L$ can be chosen so that $a$ and $b$ have odd orders. Clearly in this case $K=\langle au, bv, c\rangle$, contrary the choice of $K$. On the other hand, if $d(M)=2$ and $N=\langle u\rangle$, then there exist generators $a, b, c$ of $L$ such that $a$ is of odd order and so $K=\langle au, b, c \rangle$. This last contradiction completes the proof.\qed  

\section{Products of a periodic and a torsion-free local cyclic group} 

Recall that a group $G$ has finite torsion-free rank if it has a series of finite length whose factors are either periodic or infinite cyclic. The number $r_0(G)$ of infinite cyclic factors in such a series is an invariant of $G$ called its torsion-free rank. In this section, we describe the structure of a group $ G = AB $ with locally cyclic subgroups $ A $ and $ B $, the first of which is periodic and the other non-trivial torsion-free. Clearly $r_0(B)=1$ and we note first that  $r_0(G)=1$. 

\begin{lemma}\label{8} Let $G=AB$ be a group with subgroups $A$ and $B$ such that  $A$ is periodic abelian and $B$ is non-trivial torsion-free locally cyclic. Then $r_0(G)=1$. \end{lemma} 

\proof It was proved by D. Zaitsev in \cite{Zai_80}, Theorem 3.7 (see also \cite{AFG}, Lemma 7.1.2) that there exists a non-trivial normal subgroup of $G$ contained in $A$ or $B$. Therefore $G$ has the normal series $A_0 < A_0B_0 < G$ in which $A_0$ is the core of $A$ in $G$ and $B_0$ is the core of $B$ in $G$ modulo $A_0$. As is easily seen, the factors $A_0$ and $G/A_0B_0$ are periodic and the factor group $A_0B_0/A_0$ is isomorphic to $B_0$. Thus $r_0(G)=r_0(B)=1$, as claimed.\qed\medskip

The following lemma is a consequence of the well-known theorem of I. Schur on the finiteness of the derived subgroup of a group that is finite over its center (see \cite{Rob}, Corollary to Theorem 4.12).

\begin{lemma}\label{sh}  If a group $G$ contains a central subgroup $Z$ such that the factor group $G/Z$ is locally finite, then the derived subgroup of $G$ is locally finite.\end{lemma}

\begin{theorem}\label{9} Let the group $G=AB$ be the product of two locally cyclic subgroups $A$ and $B$ such that $A$ is periodic and $B$ is non-trivial torsion-free. Then one of the following statements hold.\begin{itemize}
\item[1)] the subgroup $A$ is normal in $G$ and so $G=A\rtimes B$;
\item[2)] $A=A_1\langle a\rangle$ with $a^2\in A_1$, the subgroup $A_1$ is normal in $G$ and $G=(A_1\rtimes B)\langle a\rangle$ with $b^a=b^{-1}\phi(b)$ for all $b\in B$, where $\phi : B \to A_1$ is a derivation of $B$ into $A_1$. \end{itemize} \end{theorem} 

\proof It is easy to see that each periodic normal subgroup $H$ of $G$  is contained in $A$, because  $AH=A(AH\cap B)$ and $AH\cap B=1$. Therefore the core $A_1=\cap_{g\in G}A^g$ of $A$ in $G$ is the maximal periodic normal subgroup of $G$.

Assume first that $A_1=1$ and let $B_1$ be the core of $B$ in $G$. Then $B_1\ne 1$ by the theorem of D. Zaitsev noted above and so the factor group $G/B_1$ is periodic, because it is the product of two periodic subgroups $AB_1/B_1$ and $B/B_1$. Moreover, since the centralizer $C_G(B_1)$ of $B_1$ in $G$ contains $B$, the group $G$ induces on $B_1$ a periodic group of automorphisms which is isomorphic to the factor group  $A/C_A(B_1)$. As is well known, a periodic group of automorphisms of any torsion-free locally cyclic group is of order $2$. Therefore the order of $A/C_A(B_1)$ does not exceed $2$ and hence  either $A=C_A(B_1)$ or $A=C_A(B_1)\langle a\rangle$ with $a\in A$ and $a^2\in C_A(B_1)$.

On the other hand, since the centralizer $C_G(B_1)=C_A(B_1)B$ is normal in $G$ and periodic over $B_1$, its derived subgroup $C_G(B_1)'$ is periodic by Lemma \ref{sh} and normal in $G$. Therefore $C_G(B_1)'\le A_1=1$ and hence $C_G(B_1)=C_A(B_1)\times B$. But then again $C_A(B_1)$ is normal in $G$ and so $C_A(B_1)=1$. Thus in the case $A_1=1$ we have either $A=1$ and $G=B$ or $A=\langle a\rangle$ with $a^2=1$ and $G=B\rtimes\langle a\rangle$ with $b^a=b^{-1}$ for all $b\in B$.

Finally, returning now to the general case, we derive that either $G=A\rtimes B$ or $G=(A_1\rtimes B)\langle a\rangle)$ with $b^a=\phi(b)b^{-1}$ for every $b\in B$ and some element $\phi(b)\in A_1$. Moreover, since $\phi(bc)(bc)^{-1}=(bc)^a =b^ac^a=(\phi(b)b^{-1})(\phi(c)c^{-1})=(\phi(b)\phi(c)^b)(bc)^{-1}$, it follows that $\phi(bc)=\phi(b)\phi(c)^b$ for any $b,c\in B$. The latter means in particular that the mapping $\phi : B \to A_1$ is a derivation of $B$ into $A_1$, as claimed.\qed

\section{Products of finitely many periodic local cyclic groups} 

A well-known theorem of B. Huppert cited in the Introduction says that every finite group of the form $G=A_1A_2...A_n$ with pairwise permuting cyclic subgroups $A_i$ for $1\le i\le n$ is supersoluble.
This result was later extended to products of pairwise permutable locally cyclic \v{C}ernikov groups by M. Tomkinson in \cite{T_86}. He proved that in this case $G=A_1A_2...A_n$ is a locally supersoluble \v{C}ernikov group. In this section we generalize this result to products of arbitrary periodic locally cyclic groups. Recall that a group is said to be  hyperabelian (respectively, hypercyclic) if it has an ascending series of normal subgroups with abelian (respectively cyclic) factors.

\begin{lemma}\label{10} Let $G=A_1A_2...A_n$ be the product of pairwise permutable periodic locally cyclic subgroups $A_i$. If the set $\pi= \cup_{i=1}^n\pi(A_i)$ is finite, $p$ is the largest prime in $\pi$, $P_i$ is the Sylow $p$-subgroup of $A_i$ and $Q_i$ is the $p$-complement to $P_i$ in $A_i$ for each $1\le i\le n$, then $G$ is a $\pi$-group, $P=P_1P_2\dots P_n$ is a normal Sylow $p$-subgroup of $G$ and $Q=Q_1Q_2\dots Q_n$ is a $p$-complement to $P$ in $G$. \end{lemma} 

\proof Since each of the $A_i$ is a subgroup of Pr\"ufer rank $1$, the group $G=A_1A_2\dots A_n$ is hyperabelian of finite Pr\"ufer rank by \cite{AS_16}, Theorem 3.1. Therefore, arguing by induction on $n$ and applying Corollary 3.2.7 of \cite{AFG} and Lemma 3.2 of \cite{AS_16}, we derive that $G$ is a $\pi$-group, $P=P_1...P_n$ is a Sylow $p$-subgroup of $G$ and $Q=Q_1Q_2\dots Q_n$ is a complement to $P$ in $G$. Moreover, taking into account that the subgroups $A_iA_j$ are locally supersoluble by Lemma \ref{3}, we conclude that ${P_i}^{A_j}\le P$ for all $i,j$ and so $P$ is a normal subgroup of $G$.\qed

\begin{lemma}\label{11} Let $G=A_1A_2...A_n$ be the product of pairwise permutable locally cyclic subgroups $A_i$. If the group $G$ is periodic and the set $\pi(G)$ is finite, then $G$ is hypercyclic. \end{lemma} 

\proof It is easy to see that every factor group of $G$ satisfies the hypothesis of the lemma. Therefore, in order to prove that $G$ is hypercyclic, it suffices to show that $G$ has a non-trivial cyclic normal subgroup. 

Let $p$ be the largest prime in $\pi(G)$, $P_i$ the Sylow $p$-subgroup of $A_i$ and $Q_i$ the $p$-complement to $P_i$ in $A_i$ for each $1\le i\le n$. Then $P=P_1P_2\dots P_n$ is a normal Sylow $p$-subgroup of $G$ and $Q=Q_1Q_2\dots Q_n$ is a $p$-complement to $P$ in $G$ by Lemma \ref{10}. It is clear that if $\pi(G)=\{p\}$, then the group $G=P_1P_2\dots P_n$ is hypercyclic by Lemma \ref{1}. Therefore, arguing by induction on $|\pi(G)|$, we may assume that $p>2$ and the subgroup $Q=Q_1Q_2\dots Q_n$ is hypercyclic.
Then every non-trivial normal subgroup of $Q$ contains a non-trivial $Q$-invariant cyclic subgroup. Since the centralizer $C_Q(P)$ of $P$ in $Q$ is a normal subgroup of $Q$, every $Q$-invariant cyclic subgroup of $C_Q(P)$ is normal in $G$. Thus, if $C_Q(P)\ne 1$, then the group $G$ contains a non-trivial cyclic normal subgroup. Therefore, to complete the proof of the lemma, it remains to consider the case $C_Q(P)=1$. 

It is clear that if in this case the subgroup $P$ is finite, then $Q$ and so the group $G$ are also finite, so that $G$ is supersoluble and thus hypercyclic by Huppert's theorem cited above. Suppose now that the subgroup $P$ is infinite. Then at least one of the subgroups $P_i$, say $P_1$, must be infinite and so quasicyclic. Since $PA_i\cap A_1A_i=(P\cap A_1A_i)A_i=(P_1P_i)A_i=P_1A_i$, we conclude that $P_1A_i=A_iP_1$ for every $i$. Furthermore, the subgroup $P_1P_i$ is abelian by Lemma \ref{1} and it is normal in $P_1A_i$, because $P_1P_i=P\cap P_1A_i$. As $A_i=P_i\times Q_i$, it follows that ${P_1}^{Q_i}=P_1$ and so the subgroup $P_1$ is normal in $P_1A_i$ for every $i$. Therefore, $P_1$ and the subgroup of order $p$ of $P_1$ are normal in $G=A_1A_2...A_n$, as claimed.\qed\medskip

{\bf Proof of Theorem 1.3.} Let $G=A_1A_2...A_n$ be the product of pairwise permutable periodic locally cyclic subgroups $A_i$. Then $G$ is a periodic group by Lemma \ref{10}. If the set $\pi(G)$ is finite, then the group $G$ is hypercyclic by Lemma \ref{11} and so locally supersoluble. In the other case the set $\pi(G)$ is infinite and thus it can be presented as a union $\pi(G)=\cup_{i=1}^\infty\pi_{i}$ of finite subsets $\pi_{i}$ such that $\pi_{i}\subset\pi_{i+1}$ for all $i\ge 1$. Let $P_{ij}$ be the Sylow $\pi_{i}$-subgroup of $A_j$ for $1\le j\le n$ and $G_i=P_{i1}P_{i2}\dots P_{in}$. Then $G_i$ is a Sylow $\pi_{i}$-subgroup of $G$ by Lemma \ref{5} which is locally supersoluble for each $i\ge 1$ by Lemma \ref{11}. Since $G=\cup_{i=1}^\infty G_i$, the group $G$ is also locally supersoluble, as claimed.\qed

\vskip 1cm %plus  1filll

\noindent Address of the authors:\\

\vbox{\halign to \hsize{\hbox to 0.45\hsize{#\hfil} & \hbox to
0.5\hsize{#\hfil}\cr

Bernhard Amberg           & Yaroslav P. Sysak \cr

Institut f\"ur Mathematik & Institute of Mathematics\cr

der Universit\"at Mainz & Ukrainian National Academy of Sciences\cr

D-55099 Mainz              & 01601 Kiev\cr

Germany                        & Ukraine \cr}

}


\begin{thebibliography}{99}

\bibitem{AFG} Amberg, B., Franciosi, S., de Giovanni, F., Products of Groups, The Clarendon Press, Oxford University Press, Oxford, 1992.

\bibitem{AS_16} Amberg, B., Sysak, Ya., Groups factorized by pairwise permutable abelian subgroups of finite rank, Adv. Group Theory Appl. 2 (2016), 13 - 24 

\bibitem{Har} Hartley, B., Splitting over the locally nilpotent residual for a class of locally finite groups, Quart. J. Math. Oxford 27 (1976), 395--400. 

\bibitem{Hu} Huppert, B., Endliche Gruppen I, Springer-Verlag, Berlin, 1967.

\bibitem{Hup} Huppert, B., \"Uber das Produkt von paarweise vertauschbaren zyklischen Gruppen, Math. Zeitschr. 58 (1953), 243--264.

\bibitem{Ja} Janko, Z., Finite $2$-groups with exactly one nonmetacyclic maximal subgroup,  Israel J. Math, 166 (2008), 313--347.

\bibitem{Ko} Kovacs, L. G., On finite soluble groups, Math. Zeitschr. 103 (1968), 37--39.

\bibitem{Rob} Robinson, D. J. S., A course in the theory of groups, Second Edition. Springer-Verlag, New York, 1996. 

\bibitem{S_86} Sysak, Ya., Products of locally cyclic, torsion-free groups,  Algebra i Logika  25 (1986),  672--686.  

\bibitem{T_86} Tomkinson, M. J., Products of abelian subgroups,  Arch. Math.  47 (1986), 107--112.  

\bibitem{Zai_80} Zaitsev, D. I., Products of abelian groups, Algebra i Logika 19       
 (1980), 150--172.

\end{thebibliography}
\end{document}